\documentclass[11pt]{article}
\usepackage{amsfonts}
\usepackage{mathrsfs}
\usepackage{amsmath}
\usepackage{amssymb}
\usepackage[mathscr]{eucal}
\usepackage{graphicx}
\usepackage{indentfirst}
\renewcommand{\paragraph}{\roman{paragraph}}

\def \ov{\overline}

\newtheorem{theorem}{\scshape \mdseries  Theorem}[section]
\newtheorem{lemma}[theorem]{\scshape \mdseries  Lemma}

\begin{document}

\title{\sf Wiener Index, Hyper-wiener Index, Harary Index and Hamiltonicity of graphs}
\author{Guidong Yu\thanks{Email: guidongy@163.com.
Supported by the Natural Science Foundation of Anhui Province (No.
1808085MA04), and the Natural Science Foundation of
Department of Education of Anhui Province (No. KJ2017A362).}, Lifang Ren, Gaixiang Cai\\
  {\small  \it School of Mathematics \& Computation Sciences, Anqing Normal University, Anqing 246133,
  China}}

\date{}
\maketitle

\noindent {\bf Abstract:} In this paper, we discuss the
Hamiltonicity of graphs in terms of Wiener index, hyper-Wiener index
and Harary index of their quasi-complement or complement. Firstly,
we give some sufficient conditions for an balanced bipartite graph
with given the minimum degree to be traceable and Hamiltonian,
respectively. Secondly, we present some sufficient conditions for a
nearly balanced bipartite graph with given the minimum degree to be
traceable. Thirdly, we establish some conditions for a graph with
given the minimum degree to be traceable and Hamiltonian,
respectively. Finally, we provide some conditions for a
$k$-connected graph to be Hamilton-connected and traceable for every
vertex, respectively.


\noindent {\bf MR Subject Classifications:}

\section{Introduction}
Let $G$ be a simple graph of order $n$ with vertex set
$V(G)=\{v_1,v_2,\ldots,v_n\}$ and edge set $E(G),$ denoted by
$e(G)= |E(G)|$. The distance between two vertices $v_{i}$ and
$v_{j}$ of $G$, denoted by $d_{G}(v_{i},v_{j})$, is defined as the
minimum length of the paths between $v_{i}$ and $v_{j}$ in $G$. Let
$N_{G}(v)$ be the set of vertices which are adjacent to $v$ in $G$.
The degree of $v$ is denoted by $d_{G}(v)=|N_{G}(v)|$, the minimum
degree of $G$ is denoted by $\delta(G)$. Let $X\subseteq V(G)$,
$G-X$ is the graph obtained from $G$ by deleting all vertices in
$X$. $G$ is called {\it $k$-connected} (for $k\in \mathbb{N}$) if
$|V(G)|>k$ and $G-X$ is connected for every set $X\subseteq V(G)$
with $|X| < k$. Let $G =(X, Y; E)$ be a bipartite graph with two
part sets $X,Y$. If $|X|=|Y|$, $G =(X, Y; E)$ is called an {\it
balanced bipartite graph}. If $|X|=|Y|+1$, $G =(X, Y; E)$ is called
a {\it nearly balanced bipartite graph}. For two disjoint graphs
$G_1$ and $G_2$, the union of $G_1$ and $G_2$, denoted by $G_1+G_2$,
is defined as $V(G_1+G_2)=V(G_1)\cup V(G_2)$ and
$E(G_1+G_2)=E(G_1)\cup E(G_2)$; and the join of $G_1$ and $G_2$,
denoted by $G_1\vee G_2$, is defined as $V(G_1\vee G_2)=V(G_1)\cup
V(G_2)$, and $E(G_1\vee G_2)=E(G_1+G_2)\cup\{xy:x\in V(G_1),y\in
V(G_2)\}$. Denote $K_{n}$ the complete graph on $n$ vertices,
$O_{n}$ the empty graph on $n$ vertices (without edges),
$K_{n,m}=O_{n}\vee O_{m}$ the complete bipartite graph with two
parts having $n,m$ vertices, respectively.

The {\it complement} of $G$ is denoted by
$\ov{G}=(V(\ov{G}),E(\ov{G}))$, where $V(\ov{G})=V(G)$,
$E(\ov{G})=\{xy:~x,y\in V(G), xy\not \in E(G)\}$.
The {\it quasi-complement} of $G =(X, Y; E)$ is denoted by
$\widehat{G} :=(X, Y; E')$, where $E' =\{xy : x\in X, y\in Y, xy\not
\in E \}$.

The {\it Wiener index} of a connected graph $G$, denote by $W(G)$,
which is introduced by Wiener \cite{Wiener} in 1947, is defined to
be the sum of distances between every pair of vertices in $G$. That
is
$$W(G)=\sum_{v_{i},v_{j}\in V(G)}d_{G}(v_{i},v_{j}).$$

We denote $D_{i}(G)=D_{G}(v_{i})=\sum\limits_{v_{j}\in
V(G)}d_{G}(v_{i},v_{j})$, then
$$W(G)=\frac{1}{2}\sum_{i=1}^{n}D_{i}(G).$$

The hyper-Wiener index, as a generalization of the Wiener index, is
traditionally denoted by $WW(G)$. The hyper-Wiener index of acyclic
graphs was introduced by Milan Randi$\acute{c}$ \cite{Randi} in 1993
and extended to all connected graphs by Klein et al. \cite{Klein}.
The {\it hyper-Wiener index} of a connected graph $G$ is defined as
$$ WW(G)= \frac{1}{2}(\sum_{v_{i}, v_{j}\in V(G)}d_{G}(v_{i},
v_{j})+\sum_{v_{i}, v_{j}\in V(G)}d_{G}^{2}(v_{i}, v_{j})).$$

We denote $DD_{i}(G)=DD_{G}(v_{i})=\sum\limits_{v_{j}\in
V(G)}d_{G}^{2}(v_{i}, v_{j}),$ then
$$WW(G)=\frac{1}{4}\sum_{i=1}^{n}(D_{i}(G)+DD_{i}(G)).$$

The Harary index is also a useful topological index in chemical
graph theory and has received much attention during the past
decades. This index has been introduced in 1993 by Plav\v{s}i\'{c}
et al. \cite{6} and by Ivanciuc et al. \cite{4}, independently. For
a connected graph $G$, the {\it Harary index} of $G$, denoted by
$H(G)$, is defined as
$$H(G)=\sum_{v_{i},v_{j}\in V(G)}\frac{1}{d_{G}(v_{i},v_{j})}.$$

We denote $\widetilde{D}_{G}(v_{i})=\sum\limits_{v_{j}\in
V(G)}\frac{1}{d_{G}(v_{i},v_{j})}$, then
$$H(G)=\frac{1}{2}\sum_{i=1}^{n}\widetilde{D}_{G}(v_i).$$

A {\it Hamiltonian cycle} of the graph $G$ is a cycle which
contains all vertex of $G$. A {\it Hamiltonian path} of the graph $G$ is a path which
contains all vertex of $G$. The graph $G$ is said
to be {\it Hamiltonian} if it contains a Hamiltonian cycle, and is
said to be {\it traceable} if it contains a Hamiltonian path. If
every two vertices of $G$ are connected by a Hamiltonian path, it is
said to be {\it Hamilton-connected}. A graph $G$ is {\it traceable
from a vertex} $x$ if it has a Hamiltonian $x$-path. A graph is {\it
traceable from every vertex} if it contains a Hamilton path from
every vertex. All these concepts
 belong to Hamiltonicity of graphs. The problem of deciding whether a graph has Hamiltonicity is
one of the most difficult classical problems in graph theory.

Recently, some topological indices have been applied to this
problem. Up to now, there are some references on the Wiener index,
hyper-wiener index and Harary index conditions for a graph to be
traceable, Hamiltonian, Hamilton-connected, traceable from every
vertex. We refer readers to see
\cite{2,14,1,3,Klein,5,7,15,li,8,21}. Among them, Hua and Ning
\cite{3} give conditions for an balanced bipartite graph to be
Hamiltonian in terms of Wiener index and Harary index. Cai et al.
\cite{14}, in terms of Hyper-Wiener index, give conditions for an
balanced bipartite graph to be traceable and Hamiltonian, and a
$k$-connected graphs to be Hamiltonian, respectively. Li
\cite{15,li} gives some conditions for a $k$-connected graph to be
Hamiltonian in terms of Wiener index and Harary index, respectively.
Yu et al. \cite{21} give conditions for a $k$-connected graphs to be
Hamilton-connected and traceable for every vertex and a nearly
balance bipartite graph to be traceable in terms of Wiener index,
Harary index and Hyper-Wiener index, respectively. Especially, Liu
et al. \cite{5} \cite{7} give sufficient conditions for a graph to
be traceable and Hamiltonian in terms of the Wiener index and Harary
index of its complement, and an balanced bipartite graph to be
traceable and Hamiltonian in terms of its Wiener index and Harary
index of its quasi-complement, respectively. As a continuance of the
these results, we also study the similar problems.

In this paper, we discuss the Hamiltonicity of graphs in terms of
Wiener index, hyper-wiener index and Harary index of their
quasi-complement or complement. In section 2, we present some
notations and some lemmas needed in the following. In sections 3-4,
we present some conditions for an balanced bipartite graph with
given the minimum degree to be traceable and Hamiltion,
respectively. In section 5, we give sufficient conditions for a
nearly balanced bipartite graph with given the minimum degree to be
traceable. In sections 6-7, we present some conditions for a graph
with given the minimum degree to be traceable and Hamiltonian,
respectively. In sections 8-9, we provide some conditions for a
$k$-connected graph to be Hamilton-connected and traceable for every
vertex, respectively.

\section{Preliminarie}

In this paper, when we mention a bipartite graph, we always fix its
partite sets, e.g., $O_{n,m}$ and $O_{m,n}$ are considered as
different bipartite graphs, unless $m=n$.

Let $G_{1},G_{2}$ be two bipartite graphs, with the bipartition
${X_{1},Y_{1}}$ and ${X_{2},Y_{2}}$, respectively. We use
$G_{1}\sqcup G_{2}$ to denote the graph obtained from $G_{1} +G_{2}$
by adding all possible edges between $X_{1}$ and $Y_{2}$ and
possible edges between $Y_{1}$ and $X_{2}$. We define some classes
of graphs as follows:
\begin{eqnarray*}
B_n^k&=&O_{k,n-k}\sqcup K_{n-k,k}                (1\leq k\leq n/2),\\
C_n^k&=&O_{k,n-k}\sqcup K_{n-k-1,k}              (1\leq k\leq n/2),\\
R_{n}^{k}&=&K_{k,k}+K_{n-k,n-k}        (1\leqslant k\leqslant n/2),\\
Q_n^k&=&O_{k+1,n-k}\sqcup K_{n-k-1,k}   (1\leqslant k\leqslant n/2),\\
L_n^k&=&K_1\vee(K_k+K_{n-k-1}),    (1\leqslant k\leqslant (n-1)/2),\\
N_n^k&=&K_k\vee(K_{n-2k}+kK_1),       (1\leqslant k\leqslant (n-1)/2),\\
\underline{L}_n^k&=&K_{k+1}+K_{n-k-1}, (1\leqslant k\leqslant (n-2)/2),\\
\underline{N}_n^k&=&K_k\vee(K_{n-2k-1}+(k+1)K_1),  (1\leqslant k\leqslant (n-2)/2).\\
\end{eqnarray*}
Note that $e(C_{n}^{k})=n(n-k-1)+k^{2}$ and $C_{n}^{k}$ is not traceable.

\begin{lemma} If $\widehat{G}$ be a connected balanced bipartite graph on $2n$ vertices,
then $$W(\widehat{G})\leqslant2n^3-3n^2+2n+2(n-1)e(G).$$
\end{lemma}

{\bf Proof.} Let $\widehat{G}$ be the quasi-complement of $G$. Then
\begin{eqnarray*}
W(\widehat{G})=&&\hspace{-.6cm}\sum_{u,v\in V(\widehat{G})}d_{\widehat{G}}(u,v)\nonumber\\
=&&\hspace{-.6cm}\frac{1}{2}\sum_{i=1}^{2n}D_{i}(\widehat{G})\nonumber\\
\leqslant&&\hspace{-.6cm}\frac{1}{2}(\sum_{i=1}^{2n}(d_{\widehat{G}}(x_{i})+(2n-1)(n-d_{\widehat{G}}(x_{i}))+(2n-2)(n-1)))\nonumber\\
=&&\hspace{-.6cm}4n^3-5n^2+2n-(n-1)\sum_{i=1}^{2n}d_{\widehat{G}}(x_{i})\nonumber\\
=&&\hspace{-.6cm}4n^3-5n^2+2n-(n-1)\sum_{i=1}^{2n}(n-d_{G}(x_{i}))\nonumber\\
=&&\hspace{-.6cm}2n^3-3n^2+2n+(n-1)\sum_{i=1}^{2n}d_{G}(x_{i})\nonumber\\
=&&\hspace{-.6cm}2n^3-3n^2+2n+2(n-1)e(G).\nonumber
\end{eqnarray*}
This completes the proof.\hfill
$\blacksquare$

\begin{lemma}Let $\widehat{G}$ be a connected balanced bipartite graph on $2n$ vertices,
then $$WW(\widehat{G})\leqslant2n^4-5n^3+5n^2-n+(2n^2-n-1)e(G).$$
\end{lemma}

{\bf Proof.} Let $\widehat{G}$ be the quasi-complement of $G$. Then
\begin{eqnarray*}
WW(\widehat{G})=&&\hspace{-.6cm}\frac{1}{2}\sum_{u,v\in V(\widehat{G})}(d_{\widehat{G}}(u,v)+d_{\widehat{G}}^{2}(u,v))\nonumber\\
=&&\hspace{-.6cm}\frac{1}{4}\sum_{i=1}^{2n}(D_{i}(\widehat{G})+DD_{i}(\widehat{G}))\nonumber\\
\leqslant&&\hspace{-.6cm}\frac{1}{4}\sum_{i=1}^{2n}(d_{\widehat{G}}(x_{i})+(2n-1)(n-d_{\widehat{G}}(x_{i}))+(2n-2)(n-1))\nonumber\\
&&\hspace{-.6cm}+\frac{1}{4}\sum_{i=1}^{2n}(d_{\widehat{G}}(x_{i})+(2n-1)^{2}(n-d_{\widehat{G}}(x_{i}))+(2n-2)^{2}(n-1))\nonumber\\
=&&\hspace{-.6cm}4n^4-6n^3+4n^2-n-(n^2-\frac{1}{2}n-\frac{1}{2})\sum_{i=1}^{2n}d_{\widehat{G}}(x_{i})\nonumber\\
=&&\hspace{-.6cm}4n^4-6n^3+4n^2-n-(n^2-\frac{1}{2}n-\frac{1}{2})\sum_{i=1}^{2n}(n-d_{G}(x_{i}))\nonumber\\
=&&\hspace{-.6cm}2n^4-5n^3+5n^2-n+(2n^2-n-1)e(G).\nonumber
\end{eqnarray*}
This completes the proof.\hfill
$\blacksquare$

\begin{lemma}Let $\widehat{G}$ be a connected balanced bipartite graph on $2n$ vertices,
then
$$H(\widehat{G})\geqslant\frac{4n^3-n}{4n-2}-\frac{2n-2}{2n-1}e(G).$$
\end{lemma}

{\bf Proof.} Let $\widehat{G}$ be the quasi-complement of $G$. Then
\begin{eqnarray*}
H(\widehat{G})=&&\hspace{-.6cm}\sum_{u,v\in V(\widehat{G})}\frac{1}{d_{\widehat{G}}(u,v)}=\frac{1}{2}\sum_{i=1}^{2n}\widetilde{D}_{\widehat{G}}(v_i)\nonumber\\
\geqslant&&\hspace{-.6cm}\frac{1}{2}\sum_{i=1}^{2n}(d_{\widehat{G}}(x_{i})+\frac{1}{2n-1}(n-d_{\widehat{G}}(x_{i}))+\frac{n-1}{(2n-2)})\nonumber\\
=&&\hspace{-.6cm}\frac{4n^2-n}{4n-2}+\frac{n-1}{2n-1}\sum_{i=1}^{2n}d_{\widehat{G}}(x_{i})\nonumber\\
\end{eqnarray*}
\begin{eqnarray*}
=&&\hspace{-.6cm}\frac{4n^2-n}{4n-2}+\frac{n-1}{2n-1}\sum_{i=1}^{2n}(n-d_{G}(x_{i}))\nonumber\\
=&&\hspace{-.6cm}\frac{4n^3-n}{4n-2}-\frac{n-1}{2n-1}\sum_{i=1}^{2n}d_{G}(x_{i})\nonumber\\
=&&\hspace{-.6cm}\frac{4n^3-n}{4n-2}-\frac{2n-2}{2n-1}e(G).\nonumber
\end{eqnarray*}
This completes the proof.\hfill
$\blacksquare$

\begin{lemma} Let $\overline{G}$ be a connected graph of order $n$,
then $$W(\overline{G})\leqslant\frac{1}{2}n(n-1)+(n-2)e(G).$$
\end{lemma}

{\bf Proof.}
\begin{eqnarray*}
W(\overline{G})=&&\hspace{-.6cm}\sum_{u,v\in V(\overline{G})}d_{\overline{G}}(u,v)\nonumber\\
=&&\hspace{-.6cm}\frac{1}{2}\sum_{i=1}^{n}D_{i}(\overline{G})\nonumber\\
\leqslant&&\hspace{-.6cm}\frac{1}{2}\sum_{i=1}^{n}(d_{\overline{G}}(x_{i})+(n-1)(n-1-d_{\overline{G}}(x_{i})))\nonumber\\
=&&\hspace{-.6cm}\frac{1}{2}n(n-1)^2+\frac{1}{2}n(2-n)\sum_{i=1}^{n}d_{\overline{G}}(x_{i})\nonumber\\
=&&\hspace{-.6cm}\frac{1}{2}n(n-1)^2-\frac{1}{2}n(n-2)\sum_{i=1}^{n}(n-1-d_{G}(x_{i}))\nonumber\\
=&&\hspace{-.6cm}\frac{1}{2}n(n-1)+\frac{1}{2}(n-2)\sum_{i=1}^{n}d_{G}(x_{i})\nonumber\\
=&&\hspace{-.6cm}\frac{1}{2}n(n-1)+(n-2)e(G).\nonumber
\end{eqnarray*}
This completes the proof.\hfill
$\blacksquare$

\begin{lemma} Let $\overline{G}$ be a connected graph of order $n$, then
$$WW(\overline{G})\leqslant\frac{1}{2}n(n-1)+\frac{1}{2}(n^2-n-2)e(G).$$
\end{lemma}

{\bf Proof.}
\begin{eqnarray*}
WW(\overline{G})=&&\hspace{-.6cm}\frac{1}{2}\sum_{u,v\in
V(\overline{G})}(d_{\overline{G}}(u,v)+d_{\overline{G}}^{2}(u,v))\nonumber\\
=&&\hspace{-.6cm}\frac{1}{4}\sum_{i=1}^{n}(D_{i}(\overline{G})+DD_{i}(\overline{G}))\nonumber\\
\leqslant&&\hspace{-.6cm}\frac{1}{4}\sum_{i=1}^{n}(d_{\overline{G}}(x_{i})+(n-1)(n-1-d_{\overline{G}}(x_{i})))\nonumber\\
+&&\hspace{-.6cm}\frac{1}{4}\sum_{i=1}^{n}(d_{\overline{G}}(x_{i})+(n-1)^2(n-1-d_{\overline{G}}(x_{i})))\nonumber\\
=&&\hspace{-.6cm}\frac{1}{4}n^2(n-1)^2+\frac{1}{4}(2+n-n^2)\sum_{i=1}^{n}d_{\overline{G}}(x_{i})\nonumber\\
=&&\hspace{-.6cm}\frac{1}{4}n^2(n-1)^2+\frac{1}{4}(2+n-n^2)\sum_{i=1}^{n}(n-1-d_{G}(x_{i}))\nonumber\\
=&&\hspace{-.6cm}\frac{1}{2}n(n-1)+\frac{1}{4}(n^2-n-2)\sum_{i=1}^{n}d_G(x_i)\nonumber\\
=&&\hspace{-.6cm}\frac{1}{2}n(n-1)+\frac{1}{2}(n^2-n-2)e(G).\nonumber
\end{eqnarray*}
This completes the proof.\hfill
$\blacksquare$

\begin{lemma} Let $\overline{G}$ be a connected graph of order $n$,
then
$$H(\overline{G})\geqslant\frac{1}{2}(n^2-n)-\frac{n-2}{n-1}e(G).$$
\end{lemma}

{\bf Proof.}
\begin{eqnarray*}
H(\overline{G})=&&\hspace{-.6cm}\sum_{u,v\in V(\overline{G})}\frac{1}{d_{\overline{G}}(u,v)}\nonumber\\
=&&\hspace{-.6cm}\frac{1}{2}\sum_{i=1}^{n}\widetilde{D}_{\overline{G}}(v_i)\nonumber\\
\geqslant&&\hspace{-.6cm}\frac{1}{2}\sum_{i=1}^{n}(d_{\overline{G}}(x_{i})+\frac{1}{n-1}(n-1-d_{\overline{G}}(x_{i})))\nonumber\\
=&&\hspace{-.6cm}\frac{1}{2}n+\frac{n-2}{2(n-1)}\sum_{i=1}^{n}d_{\overline{G}}(x_{i})\nonumber\\
=&&\hspace{-.6cm}\frac{1}{2}(n^2-n)-\frac{n-2}{2(n-1)}\sum_{i=1}^nd_{G}(x_{i})\nonumber\\
=&&\hspace{-.6cm}\frac{1}{2}(n^2-n)-\frac{n-2}{n-1}e(G).\nonumber
\end{eqnarray*}
This completes the proof.\hfill
$\blacksquare$

\section{Traceable of  balanced bipartite Graphs}

\begin{lemma}\cite{20} Let $k$ be an integer and $G$ be an
 balanced bipartite graph on $2n$ vertices. If $\delta(G)\geqslant k\geqslant1$,
$n\geqslant2k+3$ and $$e(G)>n(n-k-2)+(k+2)^2,$$ then $G$ is
traceable unless $G\subseteq Q_{n}^{k}$ or $k=1$, $G\subseteq
R_{n}^{1}$.
\end{lemma}

\begin{theorem} Let $k$ be an integer and $\widehat{G}$ be a connected balanced bipartite graph on $2n$ vertices. If $\delta(G)\geqslant k\geqslant1$,
$n\geqslant2k+3$, and
$$W(\widehat{G})>4n^3-(2k+9)n^2+(2k^2+10k+14)n-2(k+2)^2,$$ then $G$ is
traceable unless $k=1$ and $G\subseteq R_{n}^{1}$.
\end{theorem}

{\bf Proof.} Since
$W(\widehat{G})>4n^3-(2k+9)n^2+(2k^2+10k+14)n-2(k+2)^2$, and by
lemma 2.1, we get $e(G)>n(n-k-2)+(k+2)^2$. By lemma 3.1 we obtain
that $G$ is traceable or $G\subseteq Q_{n}^{k}$ or $k=1$,
$G\subseteq R_{n}^{1}$.

If $G\subseteq Q_{n}^{k}$. Because $\widehat{G}$ is connected, and
$\delta(G)\geqslant k\geqslant1$, we get $e(G)\leq
n(n-k-1)+k(k+1)-(n-k-1)-k<n(n-k-2)+(k+2)^2$, a contradiction.

If $k=1$, $G\subseteq R_{n}^{1}$. Because $\widehat{G}$ is
connected, and $\delta(G)\geqslant 1$, we get $k=1$ and $G\subseteq
R_{n}^{1}$.

This completes the proof.\hfill $\blacksquare$

\begin{theorem} Let $k$ be an integer and $\widehat{G}$ be a connected balanced bipartite graph on $2n$ vertices.
If $\delta(G)\geqslant k\geqslant1$, $n\geqslant2k+3$, and
$$WW(\widehat{G})>4n^4-(2k+10)n^3+(2k^2+9k+14)n^2-(k^2+3k+3)n-(k+2)^2,$$
then $G$ is traceable unless $k=1$ and $G\subseteq R_{n}^{1}$.
\end{theorem}

{\bf Proof.} Since
$WW(\widehat{G})>4n^4-(2k+10)n^3+(2k^2+9k+14)n^2-(k^2+3k+3)n-(k+2)^2,$
and by lemma 2.2 we get $e(G)>n(n-k-2)+(k+2)^2.$ By lemma 3.1 we
obtain that $G$ is traceable or $G\subseteq Q_{n}^{k}$ or $k=1$,
$G\subseteq R_{n}^{1}$. By the same discussion as the proof of
theorem 3.2, we get the result.\hfill $\blacksquare$

\begin{theorem} Let $k$ be an integer and $G$ be a connected balanced bipartite graph on $2n$ vertices. If $\delta(G)\geqslant k\geqslant1$,
$n\geqslant2k+3$, and
$$H(\widehat{G})<\frac{(4k+12)n^2-(4k^2+20k-25)n+4k^2+16k+16}{4n-2},$$
then $G$ is traceable unless $k=1$ and $G\subseteq R_{n}^{1}$.
\end{theorem}

{\bf Proof.} Since
$H(\widehat{G})<\frac{(4k+12)n^2-(4k^2+20k-25)n+4k^2+16k+16}{4n-2}$,
and by lemma 2.3, so $e(G)>n(n-k-2)+(k+2)^2.$ By lemma 3.1 we obtain
that $G$ is traceable or $G\subseteq Q_{n}^{k}$ or $k=1$,
$G\subseteq R_{n}^{1}$. By the same discussion as the proof of
theorem 3.2, we get the result. \hfill $\blacksquare$

\section{Hamiltonian of  balanced bipartite Graphs}

\begin{lemma}\cite{19} Let $k$ be an integer and $G$ be a balanced bipartite graph on $2n$ vertices. If $\delta(G)\geqslant k\geqslant1$,
$n\geqslant2k+3$ and $$e(G)>n(n-k-1)+(k+1)^2,$$ then $G$ is hamiltonian unless $G\subseteq B_{n}^{k}$.
\end{lemma}

\begin{theorem} Let $k$ be an integer and $\widehat{G}$ be a connected balanced bipartite graph on $2n$ vertices. If $\delta(G)\geqslant k\geqslant1$,
$n\geqslant2k+3$, $$W(\widehat{G})>4n^3-(2k+7)n^2+(2k+2(k+1)^2+4)n-2(k+1)^2,$$
then $G$ is hamiltonian.
\end{theorem}

{\bf Proof.} Since
$W(\widehat{G})>4n^3-(2k+7)n^2+(2k+2(k+1)^2+4)n-2(k+1)^2$, and by
lemma 2.1, thus $e(G)>n(n-k-1)+(k+1)^2$. By lemma 4.1, $G$ is
hamiltonian or $G\subseteq B_{n}^{k}$.

If $G\subseteq B_{n}^{k}$. Because $\widehat{G}$ is connected, and
$\delta(G)\geqslant k\geqslant1$, we get $e(G)\leq
(n-k)^2+nk-(n-k)-k<n(n-k-1)+(k+1)^2$, a contradiction.

This completes the proof.\hfill $\blacksquare$

\begin{theorem} Let $k$ be an integer and $\widehat{G}$ be a connected balanced bipartite graph on $2n$ vertices. If $\delta(G)\geqslant k\geqslant1$,
$n\geqslant2k+3$, $$WW(\widehat{G})>4n^4-(2k+8)n^3+(2k^2+5k+7)n^2-(k^2+k+1)n-(k+1)^2,$$
then $G$ is hamiltonian.
\end{theorem}

{\bf Proof.} Since
$WW(\widehat{G})>4n^4-(2k+8)n^3+(2k^2+5k+7)n^2-(k^2+k+1)n-(k+1)^2$,
and by lemma 2.2 we get $e(G)>n(n-k-1)+(k+1)^2$. By lemma 4.1, $G$
is hamiltonian unless $G\subseteq B_{n}^{k}$. By the same discussion
as the proof of theorem 4.2, we get the result. \hfill
$\blacksquare$

\begin{theorem} Let $k$ be an integer and $\widehat{G}$ be a connected balanced bipartite graph on $2n$ vertices. If $\delta(G)\geqslant k\geqslant1$,
$n\geqslant2k+3$, $$H(\widehat{G})<\frac{(4k+8)n^2-(4k^2+12k+9)n+4k^2+8k+4}{4n-2},$$
then $G$ is hamiltonian.
\end{theorem}

{\bf Proof.} Since
$H(\widehat{G})<\frac{(4k+8)n^2-(4k^2+12k+9)n+4k^2+8k+4}{4n-2}$, and
by lemma 2.3, so $e(G)>n(n-k-1)+(k+1)^2$. By lemma 4.1, $G$ is
hamiltonian unless $G\subseteq B_{n}^{k}$. By the same discussion as
the proof of theorem 4.2, we get the result.\hfill $\blacksquare$

\section{Traceable of nearly balanced bipartite Graphs}

\begin{lemma}(Yu, Fang and Fan \cite{10}) Let $G$ be a nearly balanced bipartite graph of order $2n-1$. If
$\delta(G)\geqslant k\geqslant 1$, $n\geqslant 2k+1$ , and
$$e(G)>n(n-k-2)+(k+1)^{2},$$ then $G$ is traceable unless $G\subseteq
C_{n}^{k}$.
\end{lemma}

\begin{theorem} Let $\widehat{G}$ be a connected nearly balanced bipartite graph of order $2n-1$, where $n\geqslant2k+1,$
$\delta(G)\geqslant k\geqslant1$. If $$W(\widehat{G})>4n^3-(2k+14)n^2+(4k+2(k+1)^2+16)n-4(k+1)^2-4,$$ then
$G$ is traceable.
\end{theorem}

{\bf Proof.} Let $G=G[X,Y]$, where $X=\{x_{1},x_{2},\cdots,x_{n}\}$,
$Y=\{y_{1},y_{2},\cdots,y_{n-1}\}$.
\begin{eqnarray*}
W(\widehat{G})=&&\hspace{-.6cm}\sum_{u,v\in V(\widehat{G})}d_{\widehat{G}}(u,v)\nonumber\\
=&&\hspace{-.6cm}\frac{1}{2}\sum_{i=1}^{2n-1}D_{i}(\widehat{G})\nonumber\\
\leqslant&&\hspace{-.6cm}\frac{1}{2}(\sum_{i=1}^{n}(d_{\widehat{G}}(x_{i})+(2n-3)(n-1-d_{\widehat{G}}(x_{i}))+(2n-2)(n-1)))\nonumber\\
&&\hspace{-.6cm}+\frac{1}{2}(\sum_{i=1}^{n-1}(d_{\widehat{G}}(y_{i})+(2n-3)(n-d_{\widehat{G}}(y_{i}))+(2n-4)(n-2)))\nonumber\\
=&&\hspace{-.6cm}4n^3-12n^2+12n-4-(n-2)(\sum_{i=1}^{n}d_{\widehat{G}}(x_{i})+\sum_{i=1}^{n-1}d_{\widehat{G}}(y_{i}))\nonumber\\
=&&\hspace{-.6cm}2n^3-6n^2+8n-4+(n-2)(\sum_{i=1}^{n}d_{G}(x_{i})+\sum_{i=1}^{n-1}d_{G}(y_{i}))\nonumber\\
=&&\hspace{-.6cm}2n^3-6n^2+8n-4+2(n-2)e(G),\nonumber
\end{eqnarray*}
where $d_{\widehat{G}}(x_i)=n-1-d_{G}(x_i)$ and
$d_{\widehat{G}}(y_i)=n-d_{G}(y_i)$. Because
$W(\widehat{G})>4n^3-(2k+14)n^2+(4k+2(k+1)^2+16)n-4(k+1)^2-4$, we
get $e(G)>n(n-k-2)+(k+1)^{2}$. By lemma 5.1, we obtain that $G$ is
traceable or $G\subseteq C_{n}^{k}$.

If $G\subseteq C_{n}^{k}$. Because $\widehat{G}$ is connected, and
$\delta(G)\geqslant k\geqslant1$, we get $e(G)\leq
k^2+n(n-k-1)-(n-k-1)-k<n(n-k-2)+(k+1)^{2}$, a contradiction.

This completes the proof. \hfill $\blacksquare$

\begin{theorem} Let $\widehat{G}$ be connected a nearly balanced bipartite graph of order $2n-1$, where $n\geqslant2k+1,$
$\delta(G)\geqslant k\geqslant1$. If
$$WW(\widehat{G})>4n^4-(2k+18)n^3+(2k^2+9k+\frac{65}{2})n^2-(5k^2+12k+\frac{53}{2})n+2k^2+4k+8,$$
then $G$ is traceable.
\end{theorem}

{\bf Proof.} Let $G=G[X,Y]$, where $X=\{x_{1},x_{2},\cdots,x_{n}\}$,
$Y=\{y_{1},y_{2},\cdots,y_{n-1}\}$.
\begin{eqnarray*}
WW(\widehat{G})=&&\hspace{-.6cm}\frac{1}{2}\sum_{u,v\in
V(\widehat{G})}(d_{\widehat{G}}(u,v)+d_{\widehat{G}}^{2}(u,v))\nonumber\\
=&&\hspace{-.6cm}\frac{1}{4}\sum_{i=1}^{2n-1}(D_{i}(\widehat{G})+DD_{i}(\widehat{G}))\nonumber\\
\leqslant&&\hspace{-.6cm}\frac{1}{4}\sum_{i=1}^{n}(d_{\widehat{G}}(x_{i})+(2n-3)(n-1-d_{\widehat{G}}(x_{i}))+(2n-2)(n-1))\nonumber\\
&&\hspace{-.6cm}+\frac{1}{4}\sum_{i=1}^{n}(d_{\widehat{G}}(x_{i})+(2n-3)^{2}(n-1-d_{\widehat{G}}(x_{i}))+(2n-2))^{2}(n-1))\nonumber\\
&&\hspace{-.6cm}+\frac{1}{4}\sum_{i=1}^{n-1}(d_{\widehat{G}}(y_{i})+(2n-3)(n-d_{\widehat{G}}(y_{i}))+(2n-4)(n-2))\nonumber\\
&&\hspace{-.6cm}+\frac{1}{4}\sum_{i=1}^{n-1}(d_{\widehat{G}}(y_{i})+(2n-3)^{2}(n-d_{\widehat{G}}(y_{i}))+(2n-4))^{2}(n-2))\nonumber\\
=&&\hspace{-.6cm}4n^4-16n^3+\frac{51}{2}n^2-\frac{39}{2}n+6+\frac{1}{4}(-4n^2+10n-4)(\sum_{i=1}^{n}d_{\widehat{G}}(x_{i})+\sum_{i=1}^{n-1}d_{\widehat{G}}(y_{i}))\nonumber\\
=&&\hspace{-.6cm}2n^4-9n^3+\frac{37}{2}n^2-\frac{35}{2}n+6+\frac{1}{4}(4n^2-10n+4)(\sum_{i=1}^{n}d_{G}(x_{i})+\sum_{i=1}^{n-1}d_{G}(y_{i}))\nonumber\\
=&&\hspace{-.6cm}2n^4-9n^3+\frac{37}{2}n^2-\frac{35}{2}n+6+(2n^2-5n+2)e(G),\nonumber
\end{eqnarray*}
where $d_{\widehat{G}}(x_i)=n-1-d_{G}(x_i)$ and
$d_{\widehat{G}}(y_i)=n-d_{G}(y_i)$. Because
$WW(\widehat{G})>4n^4-(2k+18)n^3+(2k^2+9k+\frac{65}{2})n^2-(5k^2+12k+\frac{53}{2})n+2k^2+4k+8$,
we get $e(G)>n(n-k-2)+(k+1)^{2}$. By lemma 5.1, we obtain that $G$
is traceable or $G\subseteq C_{n}^{k}$. By the same discussion as
the proof of theorem 5.2, we get the result.\hfill $\blacksquare$

\begin{theorem} Let $\widehat{G}$ be a connected nearly balanced bipartite graph of order $2n-1$, where $n\geqslant2k+1,$
$\delta(G)\geqslant k\geqslant1$.
If $$H(\widehat{G})\leqslant \frac{(4k+8)n^2-(4k^2+16k+17)n+8k^2+16k+8}{4n-6},$$ then
$G$ is traceable, unless $G\subseteq C_{n}^{k}(k\leqslant6)$.
\end{theorem}

{\bf Proof.} Let $G=G[X,Y]$, where $X=\{x_{1},x_{2},\cdots,x_{n}\}$,
$Y=\{y_{1},y_{2},\cdots,y_{n-1}\}$.
\begin{eqnarray*}
H(\widehat{G})=&&\hspace{-.6cm}\sum_{u,v\in V(G)}\frac{1}{d_{\widehat{G}}(u,v)}=\frac{1}{2}\sum_{i=1}^{n}\widetilde{D}_{\widehat{G}}(v_i)\nonumber\\
\geqslant&&\hspace{-.6cm}\frac{1}{2}\sum_{i=1}^{n}(d_{\widehat{G}}(x_{i})+\frac{1}{2n-3}(n-1-d_{\widehat{G}}(x_{i}))+\frac{n-1}{(2n-2)})\nonumber\\
&&\hspace{-.6cm}+\frac{1}{2}\sum_{i=1}^{n-1}(d_{\widehat{G}}(y_{i})+\frac{1}{2n-3}(n-d_{\widehat{G}}(y_{i}))+\frac{n-2}{(2n-4)})\nonumber\\
=&&\hspace{-.6cm}\frac{4n^2-5n}{4n-6}+\frac{n-2}{2n-3}(\sum_{i=1}^{n}d_{\widehat{G}}(x_{i})+\sum_{i=1}^{n-1}d_{\widehat{G}}(y_{i}))\nonumber\\
=&&\hspace{-.6cm}\frac{4n^3-8n^2+3n}{4n-6}-\frac{n-2}{2n-3}(\sum_{i=1}^{n}d_{G}(x_{i})+\sum_{i=1}^{n-1}d_{G}(y_{i}))\nonumber\\
=&&\hspace{-.6cm}\frac{4n^3-8n^2+3n}{4n-6}-\frac{2(n-2)}{2n-3}e(G),\nonumber
\end{eqnarray*}
where $d_{\widehat{G}}(x_i)=n-1-d_{G}(x_i)$ and
$d_{\widehat{G}}(y_i)=n-d_{G}(y_i)$. Because
$$H(\widehat{G})\leqslant
\frac{(4k+8)n^2-(4k^2+16k+17)n+8k^2+16k+8}{4n-6},$$ we get
$e(G)>n(n-k-2)+(k+1)^{2}$. By lemma 5.1, we obtain that $G$ is
traceable or $G\subseteq C_{n}^{k}$.
Note $k\geqslant 7$, $H(\widehat{C_n^k})=\frac{1}{4}(n^2+2kn-n-2k^2)>\frac{(4k+8)n^3-(4k^2+24k+33)n^2+(16k^2+48k+44)n-16k^2-32k-19}{4n^2-14n+12}$.
If $G\subseteq C_n^k$, then $H(\widehat{G})>H(\widehat{C_n^k})$. We get the result. \hfill $\blacksquare$

\section{Traceable of Graphs}

\begin{lemma} \cite{19} Let $k$ be an integer and $G$ be a graph of order $n\geqslant6k+10$.
If $\delta(G)\geqslant k$ and $e(G)>\binom{n-k-2}{2}+(k+1)(k+2),$
then $G$ is traceable, unless $G\subseteq\underline{L}_{n}^{k}$ or
$\underline{N}_{n}^{k}$.
\end{lemma}

\begin{theorem} Let $k$ be an integer and $\overline{G}$ be a connected graph of order $n\geqslant6k+10$.
If $\delta(G)\geqslant k$ and
$$W(\overline{G})>\frac{1}{2}(n^3-(2k+6)n^2+(3k^2+15k+19)n-6k^2-22k-20),$$
then $G$ is traceable.
\end{theorem}

{\bf Proof.} Since
$W(\overline{G})>\frac{1}{2}(n^3-(2k+6)n^2+(3k^2+15k+19)n-6k^2-22k-20)$,
by lemma 2.4, we get $e(G)>\binom{n-k-2}{2}+(k+1)(k+2)$. By lemma
6.1, we obtain that $G$ is traceable unless
$G\subseteq\underline{L}_{n}^{k}$ or $\underline{N}_{n}^{k}$.

If $G\subseteq\underline{L}_{n}^{k}$. Note that
$W(\overline{\underline{L}_n^k})=n^2-2n-kn+k^2+2k+1$. Then if
$G\subseteq\underline{L}_{n}^{k}$, we have $W(\overline{G})\leqslant
W(\overline{\underline{L}_n^k})<\frac{1}{2}(n^3-(2k+6)n^2+(3k^2+15k+19)n-6k^2-22k-20)$,
a contradiction.

If $G\subseteq\underline{N}_{n}^{k}$. Note that
$W(\overline{\underline{N}_n^k})=\frac{1}{2}(2n^2-4n-6kn+5k^2+7k+2).$
Then if $G\subseteq\underline{N}_{n}^{k}$, we have $W(\overline{G})\leqslant W(\overline{\underline{N}_n^k})<
\frac{1}{2}(n^3-(2k+6)n^2+(3k^2+15k+19)n-6k^2-22k-20)$,
a contradiction.

This completes the proof.\hfill $\blacksquare$

\begin{theorem} Let $k$ be an integer and $\overline{G}$ be a connected graph of order $n\geqslant6k+10$.
If $\delta(G)\geqslant k$ and $$WW(\overline{G})>\frac{1}{4}n^4-(\frac{1}{2}k+\frac{3}{2})n^3+(\frac{3}{4}k^2+\frac{13}{4}k+\frac{15}{4})n^2
-(\frac{3}{4}k^2+\frac{7}{4}k+\frac{1}{2})n-\frac{3}{2}k^2-\frac{11}{2}k-5,$$
then $G$ is traceable.
\end{theorem}

{\bf Proof.} Since
$WW(\overline{G})>\frac{1}{4}n^4-(\frac{1}{2}k+\frac{3}{2})n^3+(\frac{3}{4}k^2+\frac{13}{4}k+\frac{15}{4})n^2
-(\frac{3}{4}k^2+\frac{7}{4}k+\frac{1}{2})n-\frac{3}{2}k^2-\frac{11}{2}k-5,$
 by lemma 2.5, we get $e(G)>\binom{n-k-2}{2}+(k+1)(k+2)$. By lemma 6.1, we obtain that $G$ is traceable unless $G\subseteq\underline{L}_{n}^{k}$ or
$\underline{N}_{n}^{k}$.

If $G\subseteq\underline{L}_{n}^{k}$. Note that
$WW(\overline{\underline{L}_n^k})=\frac{1}{2}(3n^2-7n-4kn+4k^2+8k+4)$.
Then if $G\subseteq\underline{L}_{n}^{k}$, we have
$WW(\overline{G})\leqslant WW(\overline{\underline{L}_n^k})<
\frac{1}{4}n^4-(\frac{1}{2}k+\frac{3}{2})n^3+(\frac{3}{4}k^2+\frac{13}{4}k+\frac{15}{4})n^2
-(\frac{3}{4}k^2+\frac{7}{4}k+\frac{1}{2})n-\frac{3}{2}k^2-\frac{11}{2}k-5$,
a contradiction.

If $G\subseteq\underline{N}_{n}^{k}$. Note that
$WW(\overline{\underline{N}_n^k})=3n^2-7n-10kn+9k^2+13k+4$.
Then if $G\subseteq\underline{N}_{n}^{k}$, we have $WW(\overline{G})\leqslant WW(\overline{\underline{N}_n^k})<\frac{1}{4}n^4-(\frac{1}{2}k+\frac{3}{2})n^3+(\frac{3}{4}k^2+\frac{13}{4}k+\frac{15}{4})n^2
-(\frac{3}{4}k^2+\frac{7}{4}k+\frac{1}{2})n-\frac{3}{2}k^2-\frac{11}{2}k-5$,
a contradiction.

This completes the proof.\hfill $\blacksquare$

\begin{theorem} Let $k$ be an integer and $\overline{G}$ be a connected graph of order $n\geqslant6k+10$.
If $\delta(G)\geqslant k$ and
$$H(\overline{G})<\frac{(2k+4)n^2-(3k^2+15k+18)n+6k^2+22k+20}{2n-2},$$
then $G$ is traceable.
\end{theorem}

{\bf Proof.} Since
$H(\overline{G})<\frac{(2k+4)n^2-(3k^2+15k+18)n+6k^2+22k+20}{2n-2},$
by lemma 2.6, we get $e(G)>\binom{n-k-2}{2}+(k+1)(k+2)$. By lemma
6.1, we obtain that $G$ is traceable unless
$G\subseteq\underline{L}_{n}^{k}$ or $\underline{N}_{n}^{k}$.

If $G\subseteq\underline{L}_{n}^{k}$. Note that
$H(\overline{\underline{L}_n^k})=\frac{1}{4}(n^2+n+2kn-2k^2-2).$
Then if $G\subseteq\underline{L}_{n}^{k}$, we have
$H(\overline{G})\geqslant H(\overline{\underline{L}_n^k})>\frac{(2k+4)n^2-(3k^2+15k+18)n+6k^2+22k+20}{2n-2}$,
a contradiction.

If $G\subseteq\underline{N}_{n}^{k}$. Note that $H(\overline{\underline{N}_n^k})=\frac{1}{4}(n^2+n+2kn-2k^2-6k-2).$
Then if $G\subseteq\underline{N}_{n}^{k}$, we have $H(\overline{G})\geqslant H(\overline{\underline{N}_n^k})>\frac{(2k+4)n^2-(3k^2+15k+18)n+6k^2+22k+20}{2n-2},$
a contradiction.

This completes the proof.\hfill $\blacksquare$

\section{Hamiltonian of Graphs}

\begin{lemma} \cite{19} Let $k$ be an integer and $G$ be a graph of order $n\geqslant6k+5$.
If $\delta(G)\geqslant k$ and $$e(G)>\binom{n-k-1}{2}+(k+1)^2.$$
then $G$ is hamiltonian, unless $G\subseteq L_{n}^{k}$ or
$N_{n}^{k}$.
\end{lemma}

\begin{theorem} Let $k$ be an integer and $\overline{G}$ be a connected graph of order $n\geqslant6k+5$.
If $\delta(G)\geqslant k$ and $$W(\overline{G})>\frac{1}{2}(n^3-(2k+4)n^2+(3k^2+11k+19)n-6k^2-14k-8),$$
then $G$ is hamiltonian.
\end{theorem}

{\bf Proof.} Since
$W(\overline{G})>\frac{1}{2}(n^3-(2k+4)n^2+(3k^2+11k+19)n-6k^2-14k-8),$
by lemma 2.4, we get $e(G)>\binom{n-k-1}{2}+(k+1)^2$. By lemma 7.1,
we obtain that $G$ is hamiltonian unless $G\subseteq L_{n}^{k}$ or
$N_{n}^{k}$.

If $G\subseteq L_{n}^{k}$. Note that $W(\overline{L_n^k})=n^2-kn-3n+k^2+k+2$.
Then if $G\subseteq{L}_{n}^{k},$ we have $W(\overline{G})\leqslant W(\overline{L_n^k})<\frac{1}{2}(n^3-(2k+4)n^2+(3k^2+11k+19)n-6k^2-14k-8),$
a contradiction.

If $G\subseteq N_{n}^{k}$. Note that $W(\overline{N_n^k})=n^2-3kn-n+\frac{5}{2}k^2+\frac{3}{2}k$.
Then if $G\subseteq{N}_{n}^{k},$ we have $W(\overline{G})\leqslant W(\overline{N_n^k})<\frac{1}{2}(n^3-(2k+4)n^2+(3k^2+11k+19)n-6k^2-14k-8)$,
a contradiction.

This completes the proof.\hfill $\blacksquare$

\begin{theorem} Let $k$ be an integer and $\overline{G}$ be a connected graph of order $n\geqslant6k+5$.
If $\delta(G)\geqslant k$ and $$WW(\overline{G})>\frac{1}{4}n^4-(\frac{1}{2}k+1)n^3+(\frac{3}{4}k^2+\frac{9}{4}k+\frac{3}{2})n^2
-(\frac{3}{4}k^2+\frac{3}{4}k+\frac{1}{4})n-(\frac{3}{2}k^2+\frac{7}{2}k+2),$$
then $G$ is hamiltonian.
\end{theorem}

{\bf Proof.} Since
$WW(\overline{G})>\frac{1}{4}n^4-(\frac{1}{2}k+1)n^3+(\frac{3}{4}k^2+\frac{9}{4}k+\frac{3}{2})n^2
-(\frac{3}{4}k^2+\frac{3}{4}k+\frac{1}{4})n-(\frac{3}{2}k^2+\frac{7}{2}k+2)$,
by lemma 2.5, we get $e(G)>\binom{n-k-1}{2}+(k+1)^2$. By lemma 7.1,
we obtain that $G$ is hamiltonian unless $G\subseteq{L}_{n}^{k}$ or
${N}_{n}^{k}$.

If $G\subseteq L_{n}^{k}$. Note that $WW(\overline{L_n^k})=3n^2-4kn-9n+4k^2+4k+6$.
Then if $G\subseteq L_{n}^{k},$  we have  $WW(\overline{G})\leqslant WW(\overline{L_n^k})<\frac{1}{4}n^4-(\frac{1}{2}k+1)n^3+(\frac{3}{4}k^2+\frac{9}{4}k+\frac{3}{2})n^2
-(\frac{3}{4}k^2+\frac{3}{4}k+\frac{1}{4})n-(\frac{3}{2}k^2+\frac{7}{2}k+2),$ a contradiction.

If $G\subseteq N_{n}^{k}$. Note that $WW(\overline{N_n^k})=3n^2-10kn-3n+9k^2+5k$. Then if $G\subseteq N_{n}^{k},$ we have $WW(\overline{G})\leqslant WW(\overline{N_n^k})<\frac{1}{4}n^4-(\frac{1}{2}k+1)n^3+(\frac{3}{4}k^2+\frac{9}{4}k+\frac{3}{2})n^2
-(\frac{3}{4}k^2+\frac{3}{4}k+\frac{1}{4})n-(\frac{3}{2}k^2+\frac{7}{2}k+2),$
a contradiction.

This completes the proof.\hfill $\blacksquare$

\begin{theorem} Let $k$ be an integer and $\overline{G}$ be a connected graph of order $n\geqslant6k+5$.
If $\delta(G)\geqslant k$ and $H(\overline{G})<\frac{(2k+2)n^2-(3k^2+11k+8)n+6k^2+14k+8}{2n-2}$,
then $G$ is hamiltonian.
\end{theorem}

{\bf Proof.} Since
$H(\overline{G})\frac{(2k+2)n^2-(3k^2+11k+8)n+6k^2+14k+8}{2n-2}$, by
lemma 2.6, we get $e(G)>\binom{n-k-1}{2}+(k+1)^2$. By lemma 7.1, we
obtain that $G$ is hamiltonian unless $G\subseteq{L}_{n}^{k}$ or
${N}_{n}^{k}$.

If $G\subseteq L_{n}^{k}$. Note that $H(\overline{L_n^k})=\frac{1}{4}(n^2-3n+2kn-2k^2-2k+2)$.
Then if $G\subseteq L_{n}^{k}$, we have $H(\overline{G})\geqslant H(\overline{L_n^k})>\frac{(2k+2)n^2-(3k^2+11k+8)n+6k^2+14k+8}{2n-2},$
a contradiction.

If $G\subseteq N_{n}^{k}$.  Note that $H(\overline{N_n^k})=\frac{1}{4}(n^2-n-3k^2+k)$,
Then if $G\subseteq N_{n}^{k}$, we have $H(\overline{G})\geqslant H(\overline{N_n^k})>\frac{(2k+2)n^2-(3k^2+11k+8)n+6k^2+14k+8}{2n-2},$
a contradiction.

This completes the proof.\hfill $\blacksquare$

\section{Hamilton-connected of Graphs}

\begin{lemma} \cite{21} Let $G$ be a $k$-connected graph of order $n$, where
$k\geqslant2$. If $e(G)>\frac{n(n-1)-k(n-k-1)}{2},$ then $G$ is
Hamilton-connected.
\end{lemma}

By lemmas 2.4, 2.5, 2.6 and 8.1, and by direct computations, we get
theorems 8.2, 8.3, and 8.4, respectively.

\begin{theorem} Let $G$ be a $k$-connected graph of order $n$, $\overline{G}$ be a connected graph, where
$k \geqslant2$. If
$$W(\overline{G})>\frac{1}{2}n^3-(\frac{1}{2}k+1)n^2+\frac{1}{2}(k^2+3k+1)n-k^2-k,$$
then $G$ is Hamilton-connected.
\end{theorem}

\begin{theorem} Let $G$ be a $k$-connected graph of order $n$, $\overline{G}$ be a connected graph, where
$k \geqslant2$. If
$$WW(\overline{G})>\frac{1}{4}n^4-(\frac{1}{4}k+\frac{1}{2})n^3
+(\frac{1}{4}k^2+\frac{1}{2}k+\frac{1}{4})n^2-(\frac{1}{4}k^2-\frac{1}{4}k)n-\frac{1}{2}k^2-\frac{1}{2}k,$$
then $G$ is Hamilton-connected.
\end{theorem}

\begin{theorem} Let $G$ be a $k$-connected graph of order $n$, $\overline{G}$ be a connected graph, where
$k \geqslant2$. If
$$H(\overline{G})<\frac{(k+1)n^2-(k^2+3k+1)n+2k^2+2k}{2(n-1)},$$ then
$G$ is Hamilton-connected.
\end{theorem}

\section{Traceable from every vertex of Graphs}

\begin{lemma}\cite{21} Let $G$ be a $k$-connected graph of order $n$, where
$k\geqslant2$. If $e(G)>\frac{n(n-1)-k(n-k)}{2},$ then $G$ is
traceable from every vertex.
\end{lemma}

By lemmas 2.4, 2.5, 2.6 and 9.1, and by a direct computation, we get
theorems 9.2, theorem 9.3, theorem 9.4, respectively.

\begin{theorem} Let $G$ be a $k$-connected graph of order $n$, $\overline{G}$ be a connected graph, where
$k \geqslant2$.  If
$$W(\overline{G})>\frac{1}{2}n^3-(\frac{1}{2}k+1)n^2+\frac{1}{2}(k^2+2k+1)n-k^2,$$
then $G$ is traceable from every vertex.
\end{theorem}

\begin{theorem} Let $G$ be a $k$-connected graph of order $n$, $\overline{G}$ be a connected graph, where
$k \geqslant2$. If
$$WW(\overline{G})>\frac{1}{4}n^4-(\frac{1}{4}k+\frac{1}{2})n^3+(\frac{1}{4}k^2+\frac{1}{4}k+\frac{1}{4})n^2
-(\frac{1}{4}k^2-\frac{1}{2}k)n-\frac{1}{2}k^2,$$
then $G$ is traceable from every vertex.
\end{theorem}

\begin{theorem} Let $G$ be a $k$-connected graph of order $n$, $\overline{G}$ be a connected graph, where
$k \geqslant2$. If
$$H(\overline{G})<\frac{(k+1)n^2+(-k^2-2k-1)n+2k^2}{2n-2},$$
then $G$ is traceable from every vertex.
\end{theorem}

\end{document}